\input amstex
\documentstyle{amsppt}
\magnification=1200

\catcode`\@=12


\pagewidth{4.8125 in}
\pageheight{7 in}
\hoffset=.2 in
\voffset=.3 in

\def\leaderfill{\leaders \hbox to 1em{\hss. \hss}\hfill}
\def\bs{\backslash}

\redefine\R{{{\bold R}}}
\redefine\C{{{\bold C}}}
\redefine\Z{{{\bold Z}}}

\redefine\L{\text{L}_\C^2}

\def\fg{{\frak g}}
\def\fz{{\frak z}}
\def\fk{{\frak k}}
\def\so{{\frak{so}}}

\define\Lap{\Delta}

\def\dim{\mathop{dim}}
\def\ker{{\mathop{ker}}}

\def\al{{\alpha}}
\def\hc{{\Cal H}}
\def\bhc{{\tilde{\Cal H}}}

\def\bt{{\bar{\tau}}}

\def\th{{\bar{H}}}

\def\fb{{\bar{f}}}

\def\bt{{\bar{\tau}}}
\def\bkt{{\bar{\tau}_K}}
\def\mr{{\R^m}}
\def\kr{{\R^k}}
\def\tkr{{\R^{2k}}}
\def\rmtk{{\R^{m+2k}}}
\def\ip{{\langle\,,\,\rangle}}
\def\gb{{G_B}}
\def\bgb{{G_B}}

\def\prmtk{{(\rmtk)'}}
\def\kprmtk{{K\bs\prmtk}}
\def\ptkr{{(\tkr)'}}
\def\kptkr{{K\bs\ptkr}}

\newcount\theoremnumber
\def\grow{\advance \theoremnumber by 1}
\def\theorem#1\par{\grow
\noindent {\bf Theorem}\quad {\the\theoremnumber} .\qquad {\sl #1}\par}

\topmatter

\title Isospectral deformations of metrics on spheres \endtitle

\rightheadtext{Isospectral Riemannian metrics}

\author Carolyn S. Gordon
\endauthor

\affil Dartmouth College
\endaffil

\address 
\flushpar Carolyn S. Gordon:  Dartmouth College,  Hanover, New Hampshire, \ 03755;
\newline carolyn.s.gordon@dartmouth.edu
\endaddress

\keywords Spectral Geometry, isospectral deformations, Laplacian
\endkeywords

\subjclass
Primary 58G25; Secondary 53C20, 22E25
\endsubjclass

\abstract We construct non-trivial continuous isospectral deformations of Riemannian metrics on the ball and on the sphere in $\R^n$ for
every $n\geq 9$.  The metrics on the sphere can be chosen arbitrarily close to the round metric; in particular, they can be chosen to be
positively curved.  The metrics on the ball are both Dirichlet and Neumann isospectral and can be chosen arbitrarily close to the flat
metric.
\endabstract

\thanks The author is partially supported by a grant from the National Science Foundation.\endthanks 

\endtopmatter


\document

\heading Introduction \endheading

    To what extent does the spectrum of the Laplacian on a Riemannian 
manifold determine the geometry of the manifold?  The spectrum is known to determine certain global geometric properties such as the
dimension, volume and total scalar curvature.  Some Riemannian manifolds, such as the round spheres in low dimensions, are known to be
uniquely determined by their spectra.  Other Riemannian manifolds are known to be infinitesimally spectrally rigid; for example negatively
curved metrics on closed manifolds cannot be continuously deformed through a family of isospectral, non-isometric metrics (see \cite{GuK}
for the two-dimensional case and \cite{CS} for the general case).  On the other hand, numerous examples of isospectral
manifolds show that the spectrum does not always uniquely determine the geometry.  The main result of this article is the following theorem
(see also the somewhat stronger statements of Corollary 3.10 and Remark 3.11):

\proclaim{Theorem} (i) In every dimension $n\geq 8$, there exist continuous isospectral deformations of Riemannian metrics on the sphere
$S^n$.  The metrics can be chosen to be positively curved; in fact they can be taken to be arbitrarily close to the round metric. 

(ii)  In every dimension $n\geq 9$, there exist continuous isospectral deformations of Riemannian metrics on the ball in $\R^n$.  The metrics can be chosen
arbitrarily close to the flat metric.
\endproclaim

The metrics in part (i) of the theorem are the ones induced on the boundary of the ball by the metrics in part (ii).   In the case of the balls, we are using the
word {\it isospectral} to mean that the metrics are both Dirichlet and Neumann isospectral; i.e., the Laplacians acting on functions with 
Dirichlet or, respectively,
Neumann boundary conditions are isospectral.

The metrics on the sphere in (i) provide the first examples of isospectral deformations of positively curved metrics, in marked contrast with the infinitesimal
rigidity of negatively curved metrics on closed manifolds.  The metrics in the theorem are, to the author's knowledge, the first examples
of isospectral deformations of metrics on balls or spheres.

Z. I. Szabo \cite{Sz2} recently gave a method for constructing pairs of isospectral metrics on balls
and spheres.   Both his technique, which involves a careful analysis of the function spaces, and his examples are distinct from ours.  

Prior to 1992, all known isospectral manifolds were locally isometric; they differed only in their global geometry.  See, for example,
\cite{BGG}, \cite{Bu}, \cite{DG},
\cite{GWW}, \cite{GW1,2}, \cite{Gt1,2}, \cite{I}, \cite{M}, \cite{Su},
\cite{V} or the expository articles \cite{Be}, \cite{Br},  \cite{G3}, or \cite{GGt}. These examples reveal various global invariants which are not spectrally
determined such as the diameter and the fundamental group.

In the past eight years, many examples of isospectral manifolds with different local geometry have
appeared.  Szabo (preprint 1992) used explicit computations to construct the first examples of isospectral manifolds with
boundary having different local geometry.  The later published version \cite{Sz1} also includes closed manifolds, including a pair of
isospectral manifolds one of which is homogeneous and the other not.  In   
\cite{G1,2}, the author constructed the first isospectral closed manifolds with different local geometry; the second article uses a method based on Riemannian
submersions to prove the isospectrality of the metrics.  This method was further developed in the series of papers
\cite{GW3} (giving isospectral deformations of manifolds with boundary), 
\cite{GGSWW} (giving isospectral deformations of closed manifolds),
\cite{Sc1} (giving isospectral deformations of simply-connected manifolds), 
\cite{GSz} (giving, for example, isospectral deformations of negatively curved manifolds with boundary), and \cite{Sc2} (giving, for example, isospectral
deformations of left-invariant Riemannian metrics on compact Lie groups).   The method produces Riemannian manifolds on which a torus of
dimension at least two acts freely by isometries.  In order to obtain the deformations in Theorem 0.1, we generalize the technique,
weakening the condition that the torus action be free.  The beautiful habilitation thesis \cite{Sc2} of Dorothee Schueth provided
inspiration.
  
The paper is organized as follows.  In $\S 1$, we explain the general technique for proving isospectrality by the use of Riemannian
submersions.  In $\S 2$, we first define the class of metrics to be considered on the sphere and ball and observe that each of the metrics
admits an isometric torus action.  We then show that the open submanifold given by the union of the principal orbits is foliated by
Riemannian manifolds isometric to those studied in the earlier papers \cite{GW3} and \cite{GGSWW}; this fact is used later in $\S 3$ to
prove that the metrics in the isospectral families are not isometric.   In
$\S 3$, we construct continuous families of isospectral, non-isometric metrics within the class discussed in $\S2$.  A reader who is
interested only in the construction of the isospectral deformations and not the proofs can restrict attention to Notation 2.1, 3.1,
and 3.2, Theorems 3.3 and 3.5, Proposition 3.9, Corollary 3.10 and Remark 3.11.

The author would like to thank Dorothee Schueth for helpful suggestions and to alert the reader to a follow-up article in preparation by
Schueth in which she gives an elegant
reformulation of the technique developed here and constructs pairs of isospectral metrics on the 6-sphere.  The author would also like to
thank the Universidad Nacional de C\'ordoba, and especially Roberto Miatello and Juan Pablo Rossetti, for their hospitality while
this paper was being written.

\medskip

\heading Section 1.  Technique for constructing isospectral manifolds with different local geometry.\endheading

\definition{1.1 Background and Notation}  Let $T$ be a torus, 
let $\pi:M\to N$ be a principal $T$-bundle and endow $M$ with a Riemannian metric so that the action of $T$ is by isometries.  
Give
$N$ the induced Riemannian metric so that
$\pi$ is a Riemannian submersion.  For $a\in M$ and
$X\in T_a(M)$, write $X=X^v+X^h$ where $X^v$ is vertical (i.e., tangent to the fiber at $a$) and $X^h$ is horizontal (i.e., orthogonal to
the fiber).  Let $H_a\in T_a(M)$ denote the mean curvature vector at $a$ of the fiber through $a$.
Since $T$ acts by isometries, we have $H_{z(a)}=z_*H_a$ for all $a\in M$ and $z\in T$.  Hence $H$ is $\pi$-related to a vector 
field $\th$ on $N$.  We
will refer to $\th$ as the {\it projected mean curvature vector field of the submersion}.

Berard-Bergery and Bourguignon \cite{BB} gave a decomposition of the Laplacian $\Lap_M$ into vertical and horizontal components 
$\Lap_M =\Lap^v+\Lap^h$.  In the case of functions $f$ on
$M$ which are constant on the fibers of the submersion, so $f=\pi^*\fb$ for some function $\fb$ on $N$, then $\Lap^v(f)=0$ and  
$$\Lap_M(f)=\Lap^h(f)=\pi^*(\Lap_N(\fb)+\th(\fb)).\eqno{(1.1)}$$

If $N$ has non-trivial boundary, then $\partial M=\pi^{-1}(\partial N)$.  Since $\pi:M\to N$ is a Riemannian submersion, 
$\pi^*:C^\infty(N)\to C^\infty(M)$
maps functions on
$N$ satisfying Neumann boundary conditions to functions on $M$ satisfying Neumann boundary conditions.  Of course, 
the same is true for Dirichlet boundary conditions.
\enddefinition

\remark{Remark}  In the situation of Notation 1.1, the space $\pi^*(C^\infty(N))$ of functions on $M$ which are constant on the 
fibers of the submersion is precisely the space of
$T$-invariant functions.  Since $T$ acts by isometries, it follows that $\pi^*(C^\infty(N))$ is invariant under $\Lap_M$.  
Except in the case when the projected mean curvature vector field $\th$ is zero, the operator $\Lap_N + \th$ on $N$ is not self-adjoint;
however, it does have a discrete spectrum  since it is similar to the restriction of $\Lap_M$ to
$\pi^*(C^\infty(N))$ via the (non-unitary) isomorphism $\pi^*:C^\infty(N)\to\pi^*(C^\infty(N))$.
\endremark

\proclaim{1.2. Theorem}  Let $T$ be a torus.  Suppose $T$ acts by isometries on two compact Riemannian manifolds $M_1$ and $M_2$ and that
the action of $T$  on the principal orbits is free.  Let $M_i'$ be the union of all principal orbits in $M_i$, so $M_i'$ is an open
submanifold of $M_i$ and a principal $T$-bundle, $i=1,2$.  For each subtorus
$K$ of
$T$ of codimension at most one, suppose that there exists a diffeomorphism $\tau_K: M_1\to M_2$ which intertwines the actions of $T$ and
which induces an isometry $\bkt$ between the induced metrics on the quotient manifolds $K\bs M_1'$ and
$K\bs M_2'$.  Assume further that the isometry $\bkt$ satisfies $\bar{\tau}_{K*}(\th_K^{(1)})=\th_K^{(2)}$, 
where $\th_K^{(i)}$ is the projected mean curvature vector
field for the submersion $M_i'\to K\bs M_i'$.  (See Notation 1.1.)  Then in the case that $M_1$ and $M_2$ are closed manifolds, 
they are isospectral.  In case $M_1$ and $M_2$ have boundary, then they are Dirichlet isospectral; under the additional assumption that
$\partial(M_i)\cap M_i'$ is dense in $\partial(M_i)$ for $i=1,2$, then the manifolds are also Neumann isospectral.
\endproclaim

This theorem was proven in \cite{GSz} in the special case that the action of $T$ on $M_1$ and $M_2$ is free so that $M_i'=M_i$.  In
that case, the hypothesis that the isometry $\bkt$ arises from a diffeomorphism $\tau_K$ between $M_1$ and $M_2$ is
unnecessary.  Earlier versions, used in the papers discussed in the introduction, required that the orbits be totally geodesic.

\demo{Proof} Let $\Lap_i$ denote the Laplacian of $M_i$, and let $\L(M_i)$ denote the space of complex-valued square-integrable functions on $M_i$.  The torus
$T$ acts on
$\L(M_i)$,
$i=1,2$, and by a Fourier decomposition for this action, we have 
$$\L(M_i)=\Sigma _{\alpha\in \hat{T}}{\Cal H}_i^\alpha$$
where $\hat{T}$ consists of all characters on $T$, i.e., all homomorphisms from the group $T$ to the unit complex numbers, and 
$${\Cal H}_i^\alpha=\{f\in \L(M_i):z f=\alpha(z)f\text{ for all } z\in T\}.$$
The space of
$C^\infty$ functions on $M_i$ decomposes into its intersection with each space ${\Cal H}_i^\alpha$.  Since the torus action on $M_i$ is by
isometries, the Laplacian leaves each such subspace of smooth functions invariant.  If $M_i$ has boundary, then the space of smooth
functions satisfying Dirichlet, respectively Neumann, boundary conditions similarly decomposes into its intersections with the spaces
${\Cal H}_i^\alpha$.  

Define an equivalence relation on $\hat{T}$ by $\alpha\equiv\beta$ if $\ker(\alpha)=\ker(\beta)$.  Let $[\alpha]$ denote the 
equivalence class of $\al$ and let $[\hat{T}]$ denote the
set of equivalence classes.  Setting 
$${\Cal H}_i^{[\alpha]}=\Sigma_{\beta\in [\al]}{\Cal H}_i^\beta ,$$
then 
$$\L(M_i)=\Sigma_{[\al]\in [\hat{T}]}{\Cal H}_i^{[\alpha]}.$$
Define $Q_{[\al]}: {\Cal H}_1^{[\alpha]}\to {\Cal H}_2^{[\alpha]}$ by 
$$Q_{[\al]}(f)=f\circ\tau_{\ker(\alpha)}^{-1}$$
where $\tau_{\ker(\alpha)}$ is the map whose existence is hypothesized in the theorem, and set
$$Q=\oplus_{[\al]\in [\hat{T}]} Q_{[\al]}.$$
We will show that $Q$ intertwines the Laplacians of $M_1$ and $M_2$.

Let $\bhc_i^{[\alpha]}$ denote the subspace of $\L(M_i')$ obtained by restriction to $M_i'$ of the elements of ${\Cal H}_i^{[\alpha]}$.

For the trivial character $\alpha =1$, we have $[1]=\{1\}$, and the space $\bhc_i^{[1]}$ consists of those functions constant on the
orbits of $T$.  By equation (1.1), the projection $\pi_i:M_i'\to T\bs M_i'$ intertwines the
restriction of
$\Lap_i$ to
$\bhc_i^{[1]}$ with the operator
$\bar{\Lap}_i+\th_T^{(i)}$ acting on
$\L(T\bs M_i')$, where $\bar{\Lap}_i$ denotes the Laplacian of $T\bs M_i'$.  The isometry
$\bt_T=\bt_{\ker(\alpha)}$ gives a unitary isomorphism between
$\L(T\bs M_1')$ and
$\L(T\bs M_2')$ which intertwines
$\bar{\Lap}_1+\th_T^{(1)}$ with
$\bar{\Lap}_2+\th_T^{(2)}$.  In case the manifolds have boundary, this unitary isomorphism carries eigenfunctions of the operator
$\bar{\Lap}_1+\th_T^{(1)}$ satisfying either Dirichlet or Neumann boundary conditions to eigenfunctions of
$\bar{\Lap}_2+\th_T^{(2)}$ satisfying the same conditions.  This unitary isomorphism pulls back to the isomorphism
$\tilde{f}\to \tilde{f}\circ\tau_T^{-1}$ from
$\bhc_1^{[1]}$ to $\bhc_2^{[1]}$ satisfying $\Lap_2(\tilde{f}\circ
\tau_T^{-1}) =\Lap_1(\tilde{f})\circ\tau_T^{-1}$.  Since $M_i'$ is dense in $M_i$, we therefore have that
$\Lap_2(f\circ \tau_T^{-1}) =\Lap_1(f)\circ\tau_T^{-1}$ for smooth functions $f\in\hc_1^1$.  Thus the map $Q_{[1]}$ defined above
intertwines the Laplacians on
$\hc_1$ and
$\hc_2$. In the case of manifolds with boundary,
Dirichlet boundary conditions are preserved since
$\tau_T$ is a diffeomorphism.  Under the additional hypothesis that $\partial(M_i)\cap M_i'$ is dense in $\partial(M_i)$ for $i=1,2$,
preservation of Neumann boundary conditions follows from the statement immediately following equation (1.1).  Thus $Q_{[1]}$ carries
eigenfunctions satisfying either boundary condition to eigenfunctions satisfying the same condition.

For non-trivial $\alpha\in\hat{T}$, the kernel  of $\alpha$ is a subtorus $K$ of $T$ of
codimension one.   The space of all functions on $M_i'$ constant on the fibers of the submersion $\pi_i:M_i'\to K\bs M_i'$ is given by
$\bhc_i^{[\alpha]}\oplus \bhc_i^{[1]}$.  Arguing as in the previous paragraph, we see that the restrictions of the Laplacians
of
$M_1$ and $M_2$ to the subspaces 
${\Cal H}_1^{[\alpha]}\oplus {\Cal H}_1^{[1]}$ and, respectively, ${\Cal H}_2^{[\alpha]}\oplus {\Cal H}_2^{[1]}$ are intertwined by the map
$f\to f\circ \tau_K^{-1}$.  Moreover, since $\tau_K$ intertwines the actions of $T$, this intertwining map carries $\hc_1^{[\alpha]}$ to
$\hc_2^{[\alpha]}$; it's restriction to this subspace is precisely
$Q_{[\al]}$.
Thus $Q$ intertwines the Laplacians of $M_1$ and $M_2$.  Preservation of the appropriate boundary conditions follows exactly as in the
previous paragraph.

\enddemo

\medskip

\heading Section 2.  Construction of the metrics on the ball and sphere\endheading

\definition{2.1 Notation}(i)  Let $T$ denote the torus $\Z^k\bs \R^k$.  The Lie algebra $\fz$ of translation invariant vector fields
on $T$ is canonically identified with 
$\R^k$.  Define a representation $\rho:T\to SO(2k)$ of $T$ on $\tkr$ by diagonally 
embedding $T$ into $SO(2k)$ as the direct product of $k$ copies of the circle $SO(2)$.  Thus $\rho_*:\fz\to\so(2k)$ is a representation
of the Lie algebra $\fz$.  For
$Z\in\fz$, let
$Z^*$ denote the fundamental vector field on $\tkr$ given by $Z^*_u=\rho_*(Z)(u)$ for $u\in\tkr$.  Let
$\fz^*=\{Z^*: Z\in\fz\}$.

Given an alternating bilinear map $B:\mr\times\mr\to\fz$, define
$g_B$ to be
the unique Riemannian metric on $\R^{m+2k}=\mr\times\tkr$ satisfying the following three conditions:

\item{(a)} the canonical projection $(\mr\times\tkr,g_B)\to\mr$ is a Riemannian submersion where the
base manifold has the standard Euclidean metric;

\item{(b)} the metric induced by $g_B$ on each fiber $\tkr$ is the standard Euclidean metric;

\item{(c)}  the horizontal space at each point $(x,u)\in\mr\times\tkr$ is given by $\hc_{(x,u)}=\{y+\frac{1}{2}B(x,y)^*_u:y\in T_x(\mr)\}$. 
(Here we are using the canonical identification of $T_x(\mr)$ with $\mr$ to define $B(x,y)$.  The coefficient $\frac{1}{2}$ is included for
later convenience.)

The action $\rho$ of the torus $T$ on $\tkr$ gives rise to an action of $T$ by isometries on $(\R^{m+2k},g_B)$ preserving each fiber. 
The action is not free.  However,
$T$ does act freely on the principal orbits; these fill an open dense subset of $\R^{m+2k}$.

(ii) Let $D$, respectively $S$, denote the closed unit ball, respectively unit sphere, in $\R^{m+2k}$ relative to the Euclidean inner
product.  We continue to denote by $g_B$ both the restriction of $g_B$ to the bounded domain $D$ and the induced Riemannian metric on $S$. 
For $z\in T$, the associated isometry of $(\R^{m+2k},g_B)$ leaves both $S$ and $D$ invariant.  Thus $T$ acts by isometries on $(D,g_B)$ and
$(S,g_B)$ with the action on the principal orbits being free.

\enddefinition

In $\S 3$ we will construct continuous families of anti-symmetric
bilinear maps
$B:\mr\times\mr\to\fz$ so that the associated Riemannian metrics $g_B$ on $S$ and on $D$ are pairwise isospectral but not isometric.  The
isospectrality proof will be an easy consequence of Theorem 1.2.  To prepare for the proof in $\S 3$ that the deformations are non-trivial,
we now give an alternate description of the metric $g_B$. 

\definition{2.2 Notation and Remarks} (i)  Denote elements of the torus $T=\Z^k\bs \R^k$ by $\bar{z}$, with $z\in\R^k$.  Given an
alternating bilinear map
$B:\mr\times\mr\to\R^k$, endow $\mr\times T$ with the structure of a two-step nilpotent Lie
group by defining the group multiplication as
$$(x,\bar{z})(x',\bar{z}')=(x+x',\bar{z}+\bar{z}'+\frac{1}{2}\overline{B(x,x')}).$$
 We will denote this Lie group by
$G_B$.  Letting $\fz=\R^k$ be the Lie algebra of $T$, then the Lie algebra $\fg_B$ of $G_B$ is given by $\fg_B=\mr +\fz$ with Lie bracket
$B:\mr\times\mr\to
\fz$. In particular,
$\fz$ is central and contains the derived algebra.  

(ii) The projection $\pi:\bgb\to\mr$ gives $\bgb$ the structure of a principal $T$-bundle.  Let $$E_B= \gb\times_T \tkr$$ be the bundle over
$\mr$ with fiber
$\tkr$ associated to this principal
$T$-bundle
via the action $\rho$ of $T$ on $\tkr$ defined in 2.1.  Elements of $E_B$ are equivalence classes $[((x,\bar{z}),u)]$, with
$(x,\bar{z})\in\bgb$ and
$u\in\tkr$, under the equivalence relation $[((x,\overline{z+z'}),u)]=[((x,\bar{z}),\bar{z}'\cdot u)]$, where $\bar{z}'\cdot u$ denotes
$\rho(\bar{z}')(u)$.  Since
$\pi:\bgb\to\mr$ is a trivial bundle,
$E_B$ is diffeomorphic to $\R^{m+2k}$.   

The left-invariant vector fields on $G_B$ define vector fields on the product $G_B\times\tkr$ which are invariant under the diagonal action
of $T$.  Thus they induce vector fields on $E_B$.  Define a Riemannian metric on
$E_B$ so that the projection
$\tilde{\pi}: E_B\to\mr$ is a Riemannian submersion, the fibers have the standard Euclidean metric, and the horizontal distribution is given
by the space of left-invariant vector fields on $\bgb$ lying in the subspace $\mr$ of $\fg_{B}$, viewed as vector fields on $E_B$.  
\enddefinition

\proclaim{2.3 Proposition} 
In the notation of 2.1 and 2.2, $(\rmtk,g_B)$ is isometric to $E_B$.
\endproclaim

\demo{Proof}  Define a bundle diffeomorphism $\tau: E_B\to \rmtk$ by
$\tau([((x,\bar{z}),u)])=(x,\bar{z}\cdot u)$.  To see that
$\tau$ is an isometry between the metric on $E_B$ defined in 2.2 and the metric $g_B$ on $\rmtk$, we need only show that
$\tau$ carries the horizontal space of $E_B$ at $[((x,\bar{z}),u)]$ to that of $(\rmtk,g_B)$ at $(x,\bar{z}\cdot u)$.  

Each element $u$ of $\fg_B=\mr +\fz$ defines two vector fields on $G_B$:  a left-invariant vector field, which we temporarily denote by the
corresponding upper case letter
$U$, and a directional derivative
$D_u$ given by ignoring the group structure on $G_B$ and viewing $G_B$ as $\mr\times T$.  For the central elements $z\in\fz$ of $\fg_B$, the
two vector fields coincide.  However, for $y\in\mr$, the expression for multiplication in
$\bgb$ given in 2.2(i) shows that the vector fields $Y$ and $D_y$ are 
related by $Y_{(x,\bar{z})}=D_{y+\frac{1}{2}B(x,y)}$.  Thus letting $\tilde{Y}$ denote the vector field
on
$E_B$ associated with the left-invariant vector field $Y$ and identifying
$y$ with an element of $T_x(\mr)$, we have $\tau_{*[((x,\bar{z}),u)]}(\tilde{Y})= y+\frac{1}{2}(B(x,y))^*_{\bar{z}\cdot u}$.   Hence
$\tau$ carries the horizontal space of $E_B$ at $[((x,\bar{z}),u)]$ to that of $(\rmtk,g_B)$ at $(x,\bar{z}\cdot u)$, so $\tau$ is an
isometry.
\enddemo

\definition{2.4 Notation}
(i) Left-invariant Riemannian  metrics on the Lie group 
$G_B$ defined in 2.2 correspond to inner products on its Lie algebra
$\fg_B$.  Given an inner product $h$ on $\fz$, extend $h$ to an inner product on $\fg_B$ by taking the orthogonal direct sum with the
standard inner product on $\mr$.  Denote by
$g_h$ the associated left-invariant Riemannian metric on $G_B$. .  The projection $\pi:\bgb\to\mr$ is a
Riemannian submersion from the metric $g_h$ to the canonical metric on $\mr$.

Recall that $\fz$ may be canonically identified with $\kr$.  Given $a=(a_1,\dots,a_k)\in(\R^+)^k$, define an inner product
$h_a$ on $\kr$ so that the standard basis $\{Z_1,\dots,Z_k\}$ of $\kr$ is orthogonal and so that $\|Z_i\|_a =a_i$, where $\|\,\|_a$ is the
norm associated with $h_a$.  We will write $g_a$ to mean
$g_{h_a}$.

(ii) Let $S_r$, respectively $U_r$, denote the geodesic sphere, respectively geodesic ball, of radius $r$ in the Euclidean space $\mr$, and
let
$N_r(B,h)=\pi^{-1}(S_r)$ and $Q_r(B,h)=\pi^{-1}(U_r)$ with the Riemannian metrics induced by $g_h$.  Thus we have Riemannian submersions
$\pi:  N_r(B,h)\to S_r$ and $\pi:  Q_r(B,h)\to U_r$.   In case $h$ is the standard inner product on $\fz=\kr$, we will write
$N_r(B)$ and
$Q_r(B)$, suppressing the name of the inner product.

\enddefinition

 The manifolds $N_r(B,h)$ and $Q_r(B,h)$ were studied in \cite{GGSWW} and \cite{GW3}, respectively.

\proclaim{2.5 Proposition}  We use Notation 2.4.  For $c\in\R^+$, the map $\mu:(G_B,g_{c^2h})\to (G_{cB},g_h)$ given by
$\mu(x,\bar{z})=(x,c\bar{z})$ is both an isometry and a Lie group isomorphism.  Moreover,
$\mu$ restricts to isometries $N_r(B,c^2h)\to N_r(cB,h)$ and $Q_r(B,c^2h)\to Q_r(cB,h)$.
\endproclaim

We omit the elementary proof.

We next
describe a foliation of the dense open subest
$S'$ (respectively,
$D'$) of the sphere (respectively, ball) on which $T$ acts freely.  

\proclaim{2.6 Proposition} 
Denote elements of $\tkr$ by $u=(u_1,\dots,u_k)$ with $u_i\in\R^2$ for each $i$.  Let $u\in\tkr$ and let $a_i=\|u_i\|$,
$i=1,\dots,k$, where $\|\,\|$ is the Euclidean norm on $\R^2$.  Then:

\item{(i)} Under the action $\rho$ of $T$ on $\tkr$ defined in 2.1, the orbit of
$T$ through
$u$ is given by $T\cdot u=\{v:\|v_i\|=a_i, i=1,\dots,k\}$ .  

\item{(ii)} If $a_i\neq 0$
for all $i$, then the $T$-saturated submanifold
$L(a):=\mr\times (T\cdot u)$ of $\rmtk$ with the metric induced by $g_B$ is isometric to $(\bgb,g_a)$ where $a=(a_1,\dots,a_k)$. If
$\|u\|^2=\Sigma_{i=1}^ka_i^2 <1$,  then the intersection $S(a):=L(a)\cap S$ with the sphere $S$ in $\rmtk$ is
isometric to
$N_r(B,h_a)$, where $r=\sqrt{1-\|u\|^2}$, and the intersection $D(a):=L(a)\cap D$ is isometric to $Q_r(B,h_a)$.  
\endproclaim

\demo{Proof}  
(i) is elementary.  For (ii), let
$\tau: E_B\to (\rmtk, g_B)$ be the isometry defined in the proof of Proposition 2.3.  Then the inverse image of $L(a)$, given by
$\tau^{-1}(L(a))=\{[((x,\bar{z}),u)]:(x,\bar{z})\in G_B\}$, with the Riemannian metric induced by the metric on $E_B$ is canonically
isometric to $(G_B,g_a)$.  
The final statement of (ii) follows easily.

\enddemo

As $a$ varies, the manifolds $S(a)$ and $D(a)$ foliate $S'$ and $D'$, respectively.  These foliations will play a key role in
the proof of non-triviality of the isospectral deformations in $\S 3$.  There we will show that any isometry between the metrics on the
sphere (respectively, ball) must induce an isometry between the metrics on a suitable leaf.  We will then appeal to \cite{GW3}
(respectively,
\cite{GGSWW}) to see that the metrics on the leaf are not isometric, thus obtaining a contradiction.

\medskip

\heading Isospectral deformations\endheading

\definition{3.1 Notation and Remarks}  Let $\fz$ be a finite-dimensional vector space and $B:\mr\times\mr\to\fz$ an alternating
bilinear map.  Let $\ip$ be the standard inner product on $\mr$ and $h$ an inner product on $\fz$.  We then obtain a linear map
$j:\fz\to\so(m)$ by the condition
$$h(B(x,y),z)=\langle j(z)x,y\rangle \eqno{(3.1)}$$ for all $x,y\in\mr$ and all $z\in\fz$.
Conversely, given a linear map $j:\fz\to\so(m)$ and an inner product $h$ on $\fz$, then equation (3.1) defines an alternating bilinear map 
$B:\mr\times\mr\to\fz$.

We consider $\fz$ as in 2.1, so that $\fz$ is canonically identified with $\kr$.  Given $j:\fz\to\so(m)$, we let
$B_j:\mr\times\mr\to\fz$ be the bilinear map associated with $j$ as in equation (3.1), taking $h$ to be the standard inner product on
$\fz=\kr$.  We will use the abbreviated notation $g_j$ for the Riemannian metric $g_{B_j}$ on $\rmtk$ defined in Notation 2.1(i).
\enddefinition

\definition{3.2. Definition}  
(i) Let $\fz$ be a vector space.  A pair $j,
j^{\prime}$ of linear maps from ${\frak z}$ to $\so(m)$ will be called {\it isospectral}, denoted $j \sim
j^{\prime}$, if for each $z \in {\frak z}$, the eigenvalue spectra, with
multiplicities, of $j(z)$ and $j^{\prime}(z)$ coincide; i.e., for each $z\in\fz$, there exists an
orthogonal linear operator
$A_z$ for which 
$$ {A_z j(z) A_z^{-1} = j^{\prime}(z)}.
$$

(ii) Let $(\fz,h)$ be an inner product space.  A pair $j,
j^{\prime}$ of linear maps from ${\frak z}$ to $\so(m)$ will be called {\it equivalent}, denoted $j \simeq
j^{\prime}$, if there exist orthogonal linear maps $A$ of $\mr$ and $C$ of $\fz$ such that 
$$ A j(z) A^{-1} = j^{\prime}(C(z))
$$
for all $z\in\fz$.
\enddefinition

\proclaim{3.3 Theorem}  Let $\fz=\kr$ and let $j,j'$ be isospectral linear maps from $\fz$ to $\so(m)$ as in 3.2.  In the
notation of 2.1 and 3.1, the Riemannian metrics $g_j$ and $g_{j'}$ on the ball $D$ in $\rmtk$ are both Dirichlet and Neumann isospectral,
and the metrics
$g_j$ and
$g_{j'}$ on the sphere $S$ in $\rmtk$ are isospectral.
\endproclaim

\proclaim{3.4 Lemma}  Let $j:\fz\to\so(m)$ be a linear map.  We use the notation of 2.1 and 3.1; in particular, $\rmtk$ fibers over
$\mr$ with fibers invariant under the action of the torus $T$.  Then relative to the metric $g_j$ on $\rmtk$ defined in 3.1, the fibers are
totally geodesic Euclidean submanifolds.  Consequently, if $K$ is any subtorus of $T$ and if $(x,u)\in\rmtk$, then the mean curvature of
the orbit
$K(x,u)=(x,K\cdot u)$ in $\rmtk$ is completely determined by the mean curvature of the torus $K\cdot u$ in $\tkr$, independently of $j$.
\endproclaim

We omit the straightforward proof of the lemma.

\demo{Proof of Theorem 3.3}  We will apply Theorem 1.2.  Let
$\prmtk$ be the union of the principal orbits for the action of $T$ on $\rmtk$ defined in 2.1.  Then $T$ acts freely on $\prmtk$ and for any
subtorus
$K$ of
$T$, the quotient
$\kprmtk$ is diffeomorphic to $\mr\times (\kptkr)$, where $\ptkr$ is the union of the principal orbits for the action
$\rho$ of $K$ on
$\tkr$ 
in 2.1 (i.e., for the restriction to $K$ of the action $\rho$ of $T$).  Let
$B=B_j$ and $B'=B_{j'}$.  For $x,y\in\mr$, the fundamental vector field $B(x,y)^*$ on $\tkr$, defined as in 2.1, induces a vector field on
$\kptkr$ which we denote by $\bar{B}(x,y)$.  We define $\bar{B}'(x,y)$ similarly.

 The metric $\bar{g}_j$ on
$\kprmtk$ induced by
$g_j$ is determined by the following three properties:  
\smallskip
\item{(a)} the projection $\kprmtk\to\mr$ is a Riemannian submersion relative to the
canonical metric on the base; 

\item{(b)} the metric on the fibers $\kptkr$ is the quotient metric induced by the canonical metric on the open
subset
$\ptkr$ of $\tkr$; and 

\item{(c)} the horizontal space at $(x,\bar{u})\in \kprmtk$, where $x\in\mr$ and $\bar{u}\in \kptkr$, is given by
$\{y+\frac{1}{2}\bar{B}(x,y)_{\bar{u}}:y\in T_x(\mr)\}$, where we use the canonical identification of $T_x(\mr)$ with $\mr$ in order to
define
$\bar{B}(x,y)$.  
\smallskip
The metric
$\bar{g}_{j'}$ induced by
$g_{j'}$ is defined analogously.

First consider the case $K=T$.  In this case, the metrics $\bar{g}_{j}$ and $\bar{g}_{j'}$ on $T\bs\prmtk$ are identical.  Let $M$ denote
either the ball $D$ or the sphere $S$.  Taking $\tau_T$ to be the identity map on $M$,  then $\tau_T$ satisfies the hypothesis of Theorem
1.2.

Next suppose that
$K$ has codimension one in
$T$. Let
$Z\in\fz$ be a non-zero vector orthogonal to the Lie algebra
$\fk$ of $K$.  By the hypothesis that $j$ and
$j'$ are isospectral maps, there exists an orthogonal transformation $A=A_Z$ of $\mr$ such that $j'(Z)=Aj(Z)A^{-1}$.  Equivalently, 
$$\langle B'(Ax,Ay),Z\rangle =\langle B(x,y),Z\rangle $$
for all $x,y\in\mr$.  
Consequently,
$$\bar{B}'(Ax,Ay)=\bar{B}(x,y) \eqno{(3.2)}$$
for all
$x,y\in\mr$.  

From equation (3.2) and the defining properties (a)-(c) of the metrics $\bar{g}_j$ and $\bar{g}_{j'}$, we conclude that the
restriction to
$\prmtk$ of the diffeomorphism
$\tau_K:\rmtk\to
\rmtk$ given by
$\tau_K((x,u))=(A(x),u)$ induces an isometry $\bar{\tau}_K:(\kprmtk,\bar{g}_j)\to (\kprmtk,\bar{g}_{j'})$.  Since $\tau_K$ restricts to the
identity on the fibers of the submersion onto $\mr$, Lemma 3.4 shows that it intertwines the mean curvature vector fields for the metrics
$g_j$ and $g_{j'}$.  Now considering the restriction of $\tau_K$ to either the ball $D$ or the sphere $S$, we see that the hypothesis of
Theorem 1.2 holds.  Thus the metrics $g_j$ and $g_{j'}$ viewed either on $D$ or on $S$ are isospectral.  

\enddemo

\proclaim{Theorem 3.5}  Let $M$ denote either the ball $D$ or the sphere $S$ in $\rmtk$. Suppose that
$j:\fz\to\so(m)$ satisfies the property that there are only finitely many orthogonal maps of $\mr$ which commute 
with all the transformations $j(Z)$, $Z\in\fz$. In the notation of 2.1, 3.1 and 3.2, if
$j':\fz\to\so(m)$ is any linear map for which the metrics
$g_j$ and
$g_{j'}$ on $M$ are isometric, then $j\simeq j'$.  
\endproclaim

The hypothesis holds for generic maps $j$.  

\proclaim{Lemma 3.6}  Let $j,j':\fz\to\so(m)$ be linear maps.   In the notation of 2.4(ii):

(a) If $Q_r(B_j)$
is isometric to
$Q_r(B_{j'})$, then $j\simeq j'$.

(b)  If $j$ satisfies the genericity condition in Theorem 3.5 and if $N_r(B_j)$ is isometric to $N_r(B_{j'})$, then $j\simeq
j'$.
\endproclaim 

Parts (a) and (b) of the lemma are proven in \cite{GW3} and \cite{GGSWW}, respectively.  The genericity condition on $j$ is stated
differently in \cite{GGSWW}, Proposition 10.  To clarify, let $\fg_j =\fg_{B_j}$ be the Lie algebra with bracket $B_j$ defined in 2.2. 
Any orthogonal linear map $\alpha$ of $\mr$ which commutes with all the transformations $j(Z)$, $Z\in\fz$, extends to an
orthogonal automorphism of
$\fg_j$ by defining
$\alpha$ to be the identity map on $\fz$.  (Here we are using the word ``orthogonal'' to mean that $\alpha$ preserves the standard inner
product on $\fg_j$.)  Conversely, every orthogonal automorphism of
$\fg_j$ which restricts to the identity on
$\fz$ must be of this form.  Thus the genericity condition in Theorem 3.5 is equivalent to the condition that there are only finitely
many orthogonal automorphisms of $\fg_j$ which restrict to the identity on $\fz$.  The latter condition is slightly weaker than the
genericity condition used in 
\cite{GGSWW} in that the word ``orthogonal'' has been inserted.  However, a glance at the arguments in \cite{GGSWW} show that only this
weaker condition is actually used.

\proclaim{Lemma 3.7}  Let $M$ denote $D$, respectively $S$, and let $j,j':\fz\to\so(m)$ be any linear maps.  Suppose $\tau:(M,g_j)\to
(M,g_{j'})$ is an isometry which carries $T$ orbits to $T$ orbits.  Let $a=(a_1,\dots,a_k)\in (\R^+)^k$ satisfy $\|a\|<1$ and $a_1=\dots =
a_k$, and let $c=\sqrt{a_i}$.  In the notation of Proposition 2.6, let $M(a)$ denote $D(a)$, respectively $S(a)$.  Then
$\tau$ leaves
$M(a)$  invariant.  Thus by Proposition 2.6(ii) and Proposition 2.5, the restriction of $\tau$ to $M(a)$ gives an isometry between
$Q_r(cB_j)$ and $Q_r(cB_j')$, respectively between $N_r(cB_j)$ and $N_r(cB_j')$.
\endproclaim

\demo{Proof of Lemma 3.7}

   For
$(x,u)\in M$ with
$u=(u_1,\dots,u_k)$, the $T$ orbit $T\cdot(x,u)=(x,T\cdot u)$ is isometric to the direct product of circles of radii
$|u_1|,\dots,|u_k|$.  Thus
$\tau$ must carry the orbit $T\cdot (x,u)$ to an orbit $T\cdot (y,v)$ such that 
$(|u_1|,\dots|u_k|)=(|v_1|,\dots|v_k|)$ up to permutation.  Since for $a$ as in the lemma, $M(a)$ is the union of all those $T$-orbits
$T\cdot(x,u)$ for which $|u_1|=\dots=|u_k|=c^2$, it follows that
$\tau$ leaves
$M(a)$ invariant.
\enddemo

\proclaim{Lemma 3.8}  Suppose that $j$ satisfies the hypothesis of Theorem 3.5.  Then $T$ is a maximal torus in the full isometry group of
$(M,g_j)$.

\endproclaim

\demo{Proof of Lemma 3.8}  The proof of the Lemma will be based on the following analogous statements for the manifolds $N_r(B_j)$ and
$Q_r(B_j)$:  (These statements also use the genericity hypothesis on $j$.)

\item{(i)} $T$ is a maximal torus in the full isometry group of $N_r(B_j)$.

\item{(ii)} $T$ is a maximal torus in the full isometry group of $Q_r(B_j)$.

(i) was proven in
\cite{GGSWW} (see the proof of Proposition 10 there).  (ii) can be seen either by a direct argument or by appealing to (i) as follows:  
Suppose
$\tau$ is an isometry of $Q_r(B_j)$ which commutes with the action of $T$.  Then $\tau$ restricts to an isometry of the boundary $N_r(B_j)$
also commuting with the action of $T$.  An isometry $\tau$ of a connected manifold $Q$ is uniquely determined by its value and differential
at a single point.  Consequently, it is also determined up to two possibilities by its restriction to any submanifold $N$ of codimension
one. Indeed, for $p\in N$, the restriction of $\tau$ to $N$ determines both $\tau(p)$ and the restriction of $\tau_{*p}$ to $T_p(N)$. 
Since
$\tau_{*p}$ is an inner product space isometry, it is moreover determined up to sign on the orthogonal one
dimensional subspace of $T_p(Q)$.  (In the special case under consideration in which $N$ is the boundary of $Q$, $\tau_{*p}$ is in
fact uniquely determined by its restriction to $N$.)  Consequently (ii) follows from (i).

  We now prove the lemma.  By Lemma 3.7, $\tau$ leaves $M(a)$ invariant and its
restriction to
$M(a)$ gives an isometry of
$Q_r(c B_j)$ (if $M=D$) or of $N_r(c B_j)$ (if $M=S$).  Moreover this isometry commutes with the action of $T$.  By 3.1, $N_r(c B_j)$ is
identical to $N_r(B_{cj})$.  Note that $cj$ also satisfies the genericity hypothesis of Theorem 3.5.  Thus by statement (i) or (ii)
above, $\tau_{|M(a)}$, is determined up to countably many possibilities modulo composition with elements of $T$.  

We now use an argument similar to the proof of (ii) above to show that $\tau$ is determined up to at most two possibilities by its
restriction to $M(a)$.  Fix a point
$p$ of
$M(a)$.  The tangent space $T_p(M(a))$ has co-dimension $k$ in $T_p(M)$.  First consider
the case that $M=D$.  Each of the fundamental vector fields $Z^*$, $Z\in\fz$, defined in 2.1 is $\tau$ invariant since $\tau$ commutes with
$T$. Hence
$\tau$ leaves $\nabla_{Z^*}Z^*$ invariant for each $Z$.  As $Z$ ranges over $\fz$, the vectors $\nabla_{Z^*}Z^*(p)$ span the orthogonal
complement of $T_p(M(a))$ in
$T_p(M)$.  Thus $\tau_{*p}$ is uniquely determined by its restriction to $T_p(M(a))$ in this case.   Next in the case that $M=S$, we argue
the same way, but now the vectors $\nabla^S_{Z^*}Z^*(p)$ span only a
$(k-1)$-dimensional subspace of
$T_p(S)$ orthogonal to $T_p(M(a))$.  Hence from the restriction of $\tau$ to $M(a)$, we know the restriction of $\tau_{*p}$ to a subspace
of co-dimension one in $T_p(S)$.  As in the proof of statement (ii) above,
$\tau_{*p}$ is determined up to sign on the orthogonal one dimensional subspace.  Thus $\tau$ is determined up to at most two choices by
its restriction to $M(a)$.  The lemma follows.

\enddemo

\demo{Proof of Theorem 3.5}  Suppose
$\tau:(M,g_j)\to (M,g_{j'})$ is an isometry.  Then it must conjugate the maximal torus $T$ in the isometry group of $(M,g_j)$ to that in
the isometry group of $(M,g_{j'})$.  We conclude that the maximal torus in the isometry group of $(M,g_{j'})$ has the same dimension as
$T$, so $T$ itself is a maximal torus. Since any two maximal tori in the isometry group are conjugate, we may assume, after composing $\tau$
with an isometry of $(M,g_{j'})$, that $\tau$ conjugates the torus $T$ in the isometry group of  $(M,g_j)$ to the torus $T$ in the isometry
group of $(M,g_{j'})$.  Thus the isometry $\tau$ carries $T$ orbits to $T$ orbits.
By Lemma 3.7, it follows that
$N_r(cB_j)$ (or
$Q_r(cB_j)$) is isometric to $N_r(cB_{j'})$ (respectively, to $Q_r(cB_{j'})$).  From \cite{GGSWW} (respectively, \cite{GW3}), we
conclude that $j\simeq j'$.

\enddemo

\proclaim{3.9 Proposition  \cite{GW3}}
Let $\fz$ be an inner product space with $\dim(\fz)=2$, and let $m$ be any positive integer other
than $1,2,3,4$, or $6$. Let $W_m$ be the real vector space consisting of
all linear maps from $\fz$ to $\so(m)$. Then there is a Zariski open
subset ${\Cal O}_m$ of $W_m$ (i.e., ${\Cal O}_m$ is the complement of the zero locus 
of some non-zero polynomial
function on $W$) such that each
$j \in {\Cal O}_m$ belongs to a $d$-parameter family of isospectral,
inequivalent elements of $W_m$. Here $d\geq \frac{m(m-1)}{2} - [\frac{m}{2}]([\frac{m}{2}]+2)>1$.
In particular, $d$ is of order at least $O(m^2)$.
\endproclaim

The families of isospectral, inequivalent $j$-maps may be parameterized so that $j$ depends smoothly on the parametrization.
Although the proposition does not give any information when $m=6$, a specific example of a smooth family of isospectral, inequivalent
$j$-maps was given in \cite{GW3} when $m=6$. 

\proclaim{3.10 Corollary}  There exist continuous families of isospectral, non-isometric 
Riemannian metrics on the ball of dimension $n$ and on the sphere of dimension $n-1$ for each $n\geq 9$.  The metrics depend smoothly
on the parameter.  Given
$\epsilon >0$, the metrics on the ball can be chosen so that the sectional curvature $K$ satisfies $|K|<\epsilon$ and the metrics on the
sphere can be chosen so that $1-\epsilon<|K|<1+\epsilon$.  In particular, there exist isospectral deformations of positively curved metrics.
\endproclaim

\demo{Proof}  One of the defining properties of the Zariski set ${\Cal O}_m$ in 
\cite{GW3} is that the elements $j$ satisfy the hypothesis of Theorem 3.5. This condition also holds for the specific example in
which $m=6$.   Thus the first statement of the corollary follows immediately from Theorems 3.3 and 3.5.  For the curvature statement,
fix $j$ and consider the family of metrics $g_{cj}$, $c\in\R^+$, on $\rmtk$.  (See Notation 2.1 and
3.1).  As $c\to 0^+$, this family of metrics
converges to the Euclidean metric on $\rmtk$, and the induced family of metrics on $S$ converges to the round metric.
Given a smooth family
$\{j_t\}$ of isospectral, inequivalent maps, then
$\{cj_t\}$ is also an isospectral family of inequivalent maps for any $c\in\R^+$.  By taking $c$ small enough, we can obtain any desired
curvature bounds for the metrics associated with $cj_t$ for all $t$ in a compact subset of the parameter space.
\enddemo

\remark{3.11 Remarks}  We have actually shown that for each $j$ in the Zariski open set ${\Cal O}_m$ in Proposition 3.9 and for each $c>0$,
the Riemannian metric $g_{cj}$ on the ball, respectively sphere, in $\R^{m+4}$ lies in a continuous $d$-parameter family ${\Cal F}(cj)$ of
isospectral,  non-isometric metrics, where $d$ is at least of order $O(m^2)$.  Moreover, as $c\to 0^+$, the metric $g_{cj}$ converges to
the flat, respectively round, metric.  
\endremark


\Refs

\medskip

\widestnumber\key{GGSWW}



\ref\key Be \by P. B\'erard \paper Vari\'et\'es Riemanniennes
isospectrales non isom\'{e}triques \jour S\'em. Bourbaki
\vol705 \issue177-178 \yr1988-89\pages127--154
\endref
 

\ref\key BB
\by L. Berard-Bergery and J. P. Bourguignon
\paper Laplacians and Riemannian submersions with totally geodesic fibres
\jour Ill. J. Math.\vol 26\yr 1982\pages 181-200
\endref


\ref\key Br\by R. Brooks
\paper Constructing isospectral manifolds
\jour Amer. Math. Monthly \vol 95 \yr 1988 \pages 823--839
\endref
 
\ref\key{BGG}
\by R. Brooks, R. Gornet and W. Gustafson
\paper Mutually isospectral Riemann surfaces \jour Adv. in Math
\vol 138\yr 1998\pages 306--322\endref

 
\ref \key Bu \by P. Buser
\paper  Isospectral Riemann surfaces
\jour Ann. Inst. Fourier (Grenoble)\vol 36 \yr 1986 \pages 167--192\endref
 
\ref \key CS
\by C. Croke and V. Sharafutdinov
\paper Spectral rigidity of a compact negatively curved manifold \jour Topology
\vol 37\yr 1998\pages 1265--1273
\endref


\ref \key DG\by D. DeTurck and C. S. Gordon
\paper Isospectral deformations II: Trace formulas, metrics, and potentials
\jour Comm. Pure Appl. Math.\vol42\yr1989\pages1067--1095
\endref


\ref\key G1\manyby C. S. Gordon
\paper Isospectral closed Riemannian manifolds which are not locally
isometric
\jour J. Differential Geom.\vol 37\pages 639--649\yr1993
\endref
 
\ref\key G2\bysame
\paper Isospectral closed Riemannian manifolds which are not locally
isometric, Part II
\inbook Contemporary Mathematics: Geometry of the Spectrum
\publ Amer. Math. Soc.\vol 173\eds R. Brooks, C. Gordon, P. Perry\yr1994
\pages121--131
\endref

\ref\key G3\bysame
\paper Survey of isospectral manifolds
\inbook Handbook of Differential Geometry, F.J. E. Dillen and L. C. A. Verstraelen, eds.\vol 1
\publ Elsevier Science B.V.\yr 2000\pages 747--778
\endref
 
\ref\key GGt
\by C. Gordon and R. Gornet
\paper Spectral geometry on nilmanifolds
\inbook Progress in Inverse Spectral Geometry
\eds S. Andersson and M. Lapidus 
\publ Birkh\"auser--Verlag
\publaddr Basel
\yr 1997
\pages 23--49
\endref
 

 

\ref\key GGSWW
\by C. Gordon, R. Gornet, D. Schueth, D. Webb, and E. N. Wilson
\paper Isospectral deformations of closed Riemannian manifolds with different scalar curvature
\jour Ann. Inst. Four., Grenoble
\vol 48\yr1998\pages 593--607
\endref

\ref\key GSz\by C. Gordon and Z. I. Szabo
\paper Isospectral deformations of negatively curved Riemannian manifolds with boundary which are not locally
isometric
\endref 


\ref\key GWW\by C. S. Gordon, D. Webb, and S. Wolpert
\paper Isospectral plane domains and surfaces via Riemannian orbifolds
\jour Invent. Math.\vol 110\pages 1--22\yr1992
\endref


\ref\key GW1\manyby C. S. Gordon and E. N. Wilson
\paper Isospectral deformations of compact solvmanifolds
\jour J. Differential Geom.\yr1984\vol 19\pages 241--256
\endref
 
\ref\key GW2\bysame
\paper The spectrum of the Laplacian on Riemannian Heisenberg manifolds
\jour Michigan Math. J.\yr1986\vol 33\pages 253--271
\endref
 
\ref\key GW3\bysame
\paper Continuous families of isospectral Riemannian metrics which are not
 locally isometric
\jour J. Differential Geom.\vol 47\yr 1997\pages 504--529
\endref


\ref\key Gt1\manyby R. Gornet
\paper A new construction of isospectral Riemannian manifolds with examples
\jour Michigan Math. J. \vol 43 \yr 1996 \pages 159--188
\endref
 
\ref\key Gt2\bysame
\paper Continuous families of Riemannian manifolds isospectral on functions
but not on 1-forms
\toappear \jour J. Geom. Anal.
\endref

\ref \key GuK
\by V. Guillemin and D. Kazhdan
\paper Some inverse spectral results for negatively curved 2-manifolds
\jour Topology \vol 19 \yr 1980 \pages 301--312
\endref

 
\ref \key I\by A. Ikeda
\paper On lens spaces which are isospectral but not isometric
\jour Ann. Sci. \'Ecole Norm. Sup. (4)\vol 13\yr1980 \pages303--315
\endref
 

\ref\key M\by J. Milnor
\paper Eigenvalues of the Laplace operator on certain manifolds
\jour Proc. Nat. Acad. Sci. U.S.A.\yr1964\vol 51\page 542
\endref
 

\ref\key Sc1
\manyby D. Schueth
\paper Continuous families of isospectral metrics on simply connected
manifolds
\jour Ann. Math.\vol 149\yr 1999\pages 287-308
\endref

\ref\key Sc2\bysame
\paper Isospectral manifolds with different local geometries
\paperinfo Habilitation thesis, University of Bonn\yr 2000
\endref
 
\ref\key Su\by T. Sunada
\paper Riemannian coverings and isospectral manifolds
\jour Ann. of Math. (2)\vol 121\yr1985\pages 169--186
\endref
 
\ref\key Sz1\manyby Z. I. Szabo
\paper Locally non-isometric yet super isospectral spaces
\jour Geom. Funct. Anal.\vol 9\yr 1999\pages 185--214
\endref

\ref\key Sz2\bysame
\paper Isospectral balls and spheres
\paperinfo preprint
\endref
 
\ref\key V\by M. F. Vign\'eras
\paper Vari\'et\'es Riemanniennes isospectrales et non isom\'etriques
\jour Ann. of Math. (2)\yr1980\vol 112\pages 21--32
\endref

\endRefs
\enddocument